\documentclass[onefignum,onetabnum,onealgnum,oneeqnum]{siamart190516}
\usepackage[papersize={8.5in,11in},margin=1.0in,footskip=0.3in]{geometry}
\usepackage{amssymb}
\usepackage{amsmath}
\usepackage{enumitem}
\usepackage{microtype}

\setlist[enumerate]{leftmargin=2.6em,label=(\roman*),topsep=0.5em,parsep=0.25em}
\setlist[itemize]{leftmargin=1.5em,topsep=0.5em,,parsep=0.25em}

\definecolor{bluey}{HTML}{1f80ad}
\hypersetup{urlcolor=bluey,linkcolor=bluey,citecolor=bluey}

\usepackage{caption}
\usepackage[font={sf,small},labelsep=period,labelfont=bf,margin=0.25in]{caption}

\definecolor{bluesection}{HTML}{175280}
\makeatletter
\renewcommand\section{\@startsection{section}{1}{.25in}%
                                   {1.3ex \@plus .5ex \@minus .2ex}%
                                   {-.5em \@plus -.1em}%
                                   {\reset@font\normalsize\bfseries\color{bluesection}}}
\renewcommand\subsection{\@startsection{subsection}{2}{.25in}%
                                     {1.3ex\@plus .5ex \@minus .2ex}%
                                     {-.5em \@plus -.1em}%
                                     {\reset@font\normalsize\bfseries\color{bluesection}}}
\makeatother

\crefname{figure}{Fig.\!}{Figs.\!}
						
\usepackage{algorithm}
\usepackage[noend]{algpseudocode}
\makeatletter
\algrenewcommand\ALG@beginalgorithmic{\footnotesize}
\renewcommand{\ALG@name}{\small Algorithm}
\makeatother

\newcommand{\dig}[1]{\texttt{#1}}
\newcommand{\seqnum}[1]{\href{https://oeis.org/#1}{\rm \underline{#1}}}

\title{On two conjectures concerning the\\ternary digits of powers of two}

\author{Robert I.~Saye\thanks{Lawrence Berkeley National Laboratory, Berkeley, California, USA (\texttt{rsaye@lbl.gov})}}
\date{\today}

\begin{document}

\maketitle

\begin{abstract}
Erd\H{o}s conjectured that 1, 4, and 256 are the only powers of two whose ternary representations consist solely of \dig{0}s and \dig{1}s. Sloane conjectured that, except for $\{2^0,2^1,2^2,2^3,2^4,2^{15}\}$, every other power of two has at least one \dig{0} in its ternary representation. In this paper, numerical results are given in strong support of these conjectures. In particular, we verify both conjectures for all $2^n$ with $n \leq 2 \cdot 3^{45} \approx 5.9 \times 10^{21}$. Our approach makes use of a simple recursive construction of numbers $2^n$ having prescribed patterns in their trailing ternary digits.
\end{abstract}



\section{Introduction}

Circa 1978, Erd\H{o}s \cite{erdos1979some} conjectured that the only powers of two which do not have a \dig{2} anywhere in their ternary representation are the numbers $2^0$, $2^2$, and $2^8$. Gupta \cite{gupta1978powers} verified this to be the case for every $2^n$ with $n \leq 4373$. Extending this bound, a numerical study of Vardi \cite{vardi1991computational} confirmed no counterexamples exist for $n \leq 2 \cdot 3^{20} \approx 7 \times 10^9$. The conjecture remains open; see the additional references and analysis of Lagarias \cite{lagarias2009ternary}; see also Dimitrov \& Howe \cite{dimitrov2021powers} who study a closely related question and prove that the only powers of two whose ternary representation contains no \dig{2} and at most twenty-five \dig{1}s are the aforementioned numbers, $2^0$, $2^2$, and $2^8$. 

Similar in spirit, Sloane \cite{sloane} conjectured that, except for the numbers $\{2^0,2^1,2^2,2^3,2^4,2^{15}\}$, every other power of two contains a \dig{0} somewhere in its ternary representation. Along the same lines, one may conjecture that, for all but finite number of cases, every power of two contains a \dig{1} somewhere in its ternary representation---however, it is straightforward to show this is essentially equivalent to the conjecture of Erd\H{o}s (the exceptional cases being replaced by $2^1$, $2^3$, and $2^9$).

One may summarize all three conjectures to say that, except for a handful of small, easily predictable cases, every power of two has every possible digit somewhere in its ternary representation. Heuristically, we anticipate this to be the case because the ternary digits of powers of two are expected to be essentially random, implying that the chances of omitting a particular digit becomes vanishingly small as the overall digit count increases. However, this is far from a proof; indeed, the conjectures represent examples of exponential Diophantine equations for which few methods of attack have been found \cite{dimitrov2021powers,lagarias2009ternary}.

In this note, numerical results are given in strong support of these conjectures. In particular, we significantly extend prior verification bounds and confirm that the ternary representation of $2^n$ contains every possible ternary digit, for all $16 \leq n \leq 2 \cdot 3^{45} \approx 5.9 \times 10^{21}$. Our approach focuses on examining the trailing ternary digits of $2^n$, which can be efficiently calculated even for massive exponents. In particular, we develop a recursive algorithm to construct numbers $2^n$ having prescribed patterns in their trailing ternary digits. For example, to find a potential counterexample to  Erd\H{o}s's conjecture, one may directly enumerate in increasing order the numbers $2^n$ whose trailing digits are some combination of \dig{0}s and \dig{1}s. We note the recursive algorithm shares some aspects with the sieving method of Gupta \cite{gupta1978powers}.

As part of our analysis, we also compute the smallest power of two which has no \dig{0} in the last $k$ digits of its ternary expansion, for $k = 1, 2, \ldots$ (and similarly for trailing digits excluding \dig{1} and \dig{2}). The results agree very well with what one may expect supposing that the ternary digits of $2^n$ are essentially rolls of a three-sided die.

\section{Notation}

It is convenient to define a shorthand notation for the purposes of examining the trailing ternary digits of a number: for integers $a,b$ and $k$ a positive integer, $a \equiv_k b$ means $a \equiv b \pmod{3^k}$. In addition, $d_k(a)$ is defined as the $k^\text{th}$ ternary digit of $a$, with $d_1(a)$ being the least significant digit: more precisely, if $a = \sum_{i = 0}^n a_i 3^i$ is the ternary representation of $a$, then $d_k(a) := a_{k-1}$. As a final piece of notation, $(\cdots)_3$ indicates the digits in the ternary expansion of a number, e.g., $2^{8} = (\texttt{100111})_3$.

\section{Method}

A simple recursive construction of numbers $2^n$, having prescribed patterns in their trailing ternary digits, is made possible via the results of the following lemma; its proof is elementary, and shares some aspects with the method of Gupta \cite{gupta1978powers}. A self-contained proof of the lemma is deferred to the appendix so as to simplify the presentation.

\begin{lemma}
\label{lemm}
For a positive integer $k$, define\footnote{In fact, $u_k = \varphi(3^k)$, where $\varphi$ is the Euler totient function; Euler's theorem implies that $a^{u_k} \equiv 1 \pmod{3^k}$ for any positive integer $a$ coprime to 3.} $u_k := 2 \cdot 3^{k-1}$. Then
\begin{enumerate}
\item $u_k$ is the smallest positive integer such that $2^{u_k} \equiv_k 1$;
\item if $i, j \in \mathbb N$ and $2^i \equiv_k 2^j$, then $i$ and $j$ differ by a multiple of $u_k$;
\item if $i,j \in \mathbb N$, the $(k+1)^\text{st}$ ternary digit of $2^{i u_k + j}$ is related to that of $2^j$ via 
\[ d_{k+1}(2^{i u_k + j}) \equiv d_{k+1}(2^j) + i \cdot d_1(2^j) \pmod{3}. \]
\end{enumerate}
\end{lemma}

We demonstrate the recursive construction process by means of a generic five-digit example. (There is nothing special about the digit count of five.) Suppose we have constructed a positive integer $j < u_5$ such that the last five ternary digits of $2^{j}$ is $abcde$, for some fixed $a, \ldots, d \in \{\dig0,\dig1,\dig2\}$ and $e \in \{\dig1,\dig2\}$.  Then, for $i \in \{0,1,2\}$, we claim the numbers $j_i := i u_5 + j$ are such that $2^{j_i}$ are the smallest possible powers of two whose last six ternary digits are $\dig0abcde$, $\dig1abcde$, and $\dig2abcde$. (The order of these six-digit combinations, as $i$ iterates from 0 to 2, depends on $e$.) To see why, note that:
\begin{itemize}
\item Applying part (i) of the lemma with $k = 5$, observe that $2^{j_i} = (2^{u_5})^i 2^j \equiv_5 2^j$, and so the last five digits are preserved.
\item For $i$ held fixed, suppose that $\ell$ is a positive integer such that $\ell < j_i$ and $2^\ell$ matches the last six digits of $2^{j_i}$. Then, by part (ii) of the lemma, $j_i - \ell$ is a positive multiple of $u_6$, but this is impossible because $j_i = i u_5 + j < 2 u_5 + u_5 = u_6$. Therefore, no such $\ell$ exists and consequently $2^{j_i}$ is the smallest possible power of two whose trailing six digits match those of $2^{j_i}$.
\item Last, by part (iii), the sixth ternary digit of $2^{j_i}$ is equal (modulo 3) to the sixth digit of $2^j$ plus 0, 1, or 2 multiples of the last digit of $2^{j}$. The latter digit is either 1 or $-1$ (modulo 3), which means that, irrespective of what the sixth digit of $2^j$ is, we shall always obtain some arrangement of $\dig0abcde$, $\dig1abcde$, and $\dig2abcde$ for the last six digits of $2^{j_i}$, as $i$ iterates over $\{0,1,2\}$.
\end{itemize}

In general, we observe that adding multiples of $u_k$ to a number $j < u_k$ allow us to explicitly construct powers of two whose last $k$ digits match those of $2^j$ and whose $(k+1)^\text{st}$ digit is controlled; moreover, the recursive approach builds powers of two in the smallest order possible. As an example application, we may then use this approach to test the conjecture of Erd\H{o}s, by starting with $2^0$ (whose least significant digit is \dig1), then generate the smallest powers of two whose trailing two digits are \dig0\dig1 and \dig1\dig1, then generate the smallest powers of two whose trailing three digits are \dig0\dig0\dig1, \dig1\dig0\dig1, \dig0\dig1\dig1, and \dig1\dig1\dig1, etc. If any of these powers of two end up containing solely \dig0s and \dig1s in their ternary representation, then a counterexample to the conjecture has been discovered (provided it is not one of the trivial cases, of course); moreover, any such counterexample must be constructable by this process. 

An algorithm implementing this strategy is given in \cref{algo}. The input is $k$, the number of so-far-constructed trailing digits, the unit $u_k$ defined by \cref{lemm} and its corresponding power of two, along with an integer $j$ and its corresponding power of two. The parameter $\chi$ specifies the digit controlling the recursive construction: if $\chi = 2$ (resp., $\chi = 0$), then \cref{algo} generates powers of two whose trailing $k$ digits contain only \dig0s and \dig1s (resp., only \dig1s and \dig2s), thereby examining the conjecture of Erd\H{o}s (resp., Sloane). In particular, for $\chi = 2$, the recursion is initiated with the first power of two having $k = 1$ valid digits, i.e., ${\mathcal G}_2(k = 1, u_k = 2, 2^{u_k} = 4, j = 0, 2^j = 1)$. Meanwhile, for $\chi = 0$, the recursion is initiated via two base cases, ${\mathcal G}_0(k = 1, u_k = 2, 2^{u_k} = 4, j = 0, 2^j = 1)$ and ${\mathcal G}_0(k = 1, u_k = 2, 2^{u_k} = 4, j = 1, 2^j = 2)$. By construction, the recursive algorithm is depth-first, with a maximum depth controlled by the user-defined parameter $K$. A straightforward calculation shows that the total number of powers of two constructed by the recursive algorithm is $\Theta(2^{K})$, and that every such power is less than $2^{u_K}$. On the other hand, the total number of powers of two less than $2^{u_K}$ is $\Theta(3^K)$. In that sense, and in the context of testing the conjectures, the recursive approach exponentially reduces the search space versus the more elementary method of simply testing every power of two in increasing order.

\begin{figure}[t]
\centering
\begin{minipage}{0.44\textwidth}
\begin{algorithm}[H]
\caption{\small ${\mathcal G}_\chi(k, u_k, 2^{u_k}, j, 2^j)$.}
\label{algo}
\begin{algorithmic}[1]
	\State Determine the first occurrence of digit $\chi$ in $2^j$. \label{algodigit}
	\If{digit $\chi$ not found and $j > 16$}
		\State \textbf{output} $j$ \textit{(nontrivial counterexample found)}
	\EndIf	
	\If{$d_k(2^j) = \chi$}
		\State \textbf{return}
	\EndIf	
	\If{$k \geq K$}
		\State \textbf{return}
	\EndIf	
	\State Compute $2^{u_{k+1}} = (2^{u_k})^3$.	
	\State Execute ${\mathcal G}_\chi\bigl(k + 1, 3 u_k, 2^{u_{k+1}}, j, 2^j\bigr)$.
	\State Execute ${\mathcal G}_\chi\bigl(k + 1, 3 u_k, 2^{u_{k+1}}, j + u_k, 2^j \cdot 2^{u_k} \bigr)$.
	\State Execute ${\mathcal G}_\chi\bigl(k + 1, 3 u_k, 2^{u_{k+1}}, j + 2 u_k, 2^j \cdot (2^{u_k})^2 \bigr)$.
	\end{algorithmic}
\end{algorithm}
\end{minipage}
\caption{Recursive generation of powers of two whose trailing $k$ ternary digits are required to satisfy particular conditions.}\vspace{-1.5em}
\end{figure}

Our implementation of \cref{algo} includes the following aspects, mainly targeting its efficient execution:
\begin{itemize}
\item Except for line \ref{algodigit}, all powers of two are computed in the cyclic group modulo $3^\kappa$ for a fixed $\kappa$. In particular, we have used a tailor-made, fixed precision integer type representing a $\kappa = 54$ digit ternary number. It is implemented as a three-digit number in base $3^{18}$, with each such digit represented by a conventional 32-bit unsigned integer (\verb|uint32_t| in C++). This approach is particularly fast at computing the cubes and multiplications in \cref{algo}.
\item On line \ref{algodigit}, we first query for the occurrence of digit $\chi$ in the fixed-precision 54-ternary digit number representing $2^j$. (Here, the ``first occurrence'' essentially means $\min \{ i : d_i(j) = \chi \}$.) Although sufficiently rare, it can happen that no such digit occurs in these 54 digits, in which case we switch over to an alternative algorithm. The alternative algorithm computes $2^j$ (via exponentiation-by-squaring) in the cycling group modulo $3^\ell$ (using a similar ternary digit implementation as above), in progressively increasing lengths $\ell$, until $\chi$ is found. In essence, this method tries to compute as few of the trailing digits of $2^j$ as possible in order to find the digit $\chi$; owing to the nature of the distribution of ternary digits of powers of two, it is usually the case that not many additional digits are required. (We note that a nontrivial counterexample to the conjectures would require $\ell$ to reach the full digit length of the ternary representation of $2^j$, however this circumstance never occurred in our computational study.)
\end{itemize}

\section{Results}

Running on a modest 64-core compute server for a few days, the computational study in this work applied a maximum recursion depth of $K = 46$. This corresponds to testing the conjectures of Erd\H{o}s and Sloane against all powers $2^n$ such that $n \leq u_{46} = 2 \cdot 3^{46 - 1} \approx 5.9 \times 10^{21}$. No counterexamples were found.

As part of this study, trailing digit count ``record breakers'' were tracked. Specifically, for $\chi \in \{0,1,2\}$, we define $\rho_\chi : \mathbb N \to \mathbb N$ such that
\[ \rho_\chi(k) = \min \{ n \in \mathbb N : \text{$2^n \geq 3^{k-1}$ and $\chi$ occurs nowhere in the last $k$ ternary digits of $2^n$} \}. \]
(In particular, the powers of two must have at least $k$ ternary digits, i.e., $2^n \geq 3^{k-1}$.) As an example, $\rho_2(100) = 710982592620911336$; the last 110 ternary digits of $2^{710982592620911336}$ are
\begin{align*}
&\bigl(\texttt{0102020002100100100110011100110101011111010101010110010} \raisebox{-1.25ex}[0pt][0pt]{$\hookleftarrow$} \\
& \qquad\qquad \texttt{1000111001000101110010101011111010001110110001110111011}\bigr)_3.
\end{align*}
As another example, $\rho_0(100) = 388128961376647359$; the last 110 digits of $2^{388128961376647359}$ are
\begin{align*}
&\bigl(\texttt{2021120020121121111112111222212121111112222122221212212} \raisebox{-1.25ex}[0pt][0pt]{$\hookleftarrow$} \\
& \qquad\qquad \texttt{1122111112221212212211111121221222222111222122221212122}\bigr)_3.
\end{align*}
It is straightforward to show that $\rho_1(k) = \rho_2(k) + 1$ for all $k$. This is because $2^n$ ends in a sequence of \dig0s and \dig2s if and only if $2^{n-1}$ ends in a sequence of \dig0s and \dig1s; moreover, the maximal number of trailing non-\dig{1} digits (for the former) and non-\dig{2} digits (for the latter) are exactly the same. As a result, we only consider $\rho_0$ and $\rho_2$ in the following analysis.

\begin{figure}[t]
\centering
\includegraphics[scale=0.9]{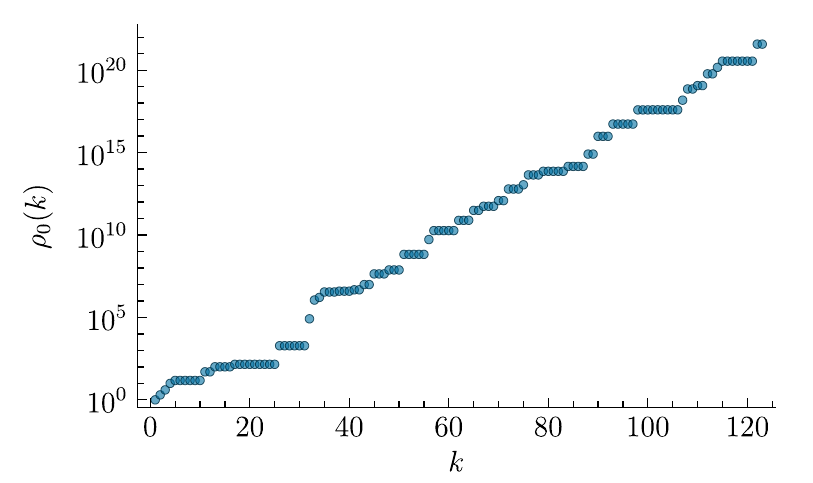}
\includegraphics[scale=0.9]{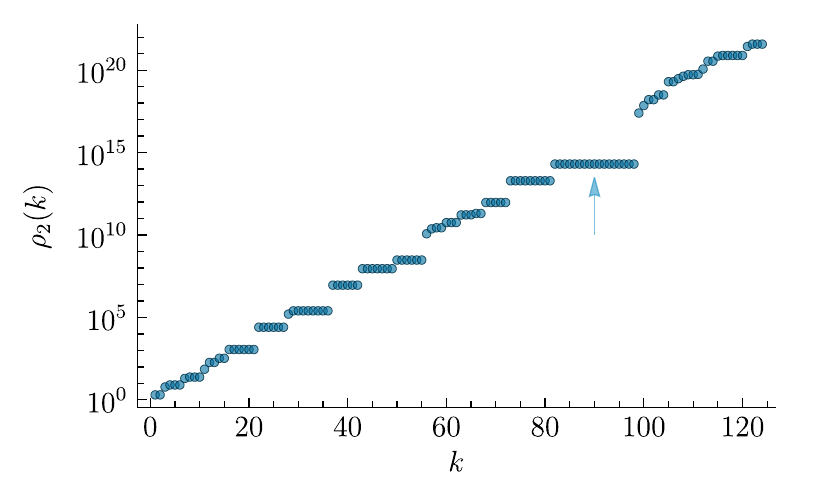}
\caption{Plots of $\rho_0$ (resp., $\rho_2$), defined as the smallest integer $n$ such that the digit \dig{0} (resp., \dig2) occurs nowhere in the last $k$ ternary digits of $2^n$. The arrow points to the instance where $\rho_2(k) = 201015414581294$ for all $82 \leq k \leq 98$.}
\label{fig1}\vspace{-1.5em}
\end{figure}

\cref{fig1} plots $\rho_0$ and $\rho_2$ as a function of $k$. We observe that $\rho_\chi(k)$ grows approximately exponentially with $k$. The longer horizontal steps correspond to the record breakers which have, roughly speaking, an uncharacteristic number of trailing non-$\chi$ digits. One notable example is $n = 201015414581294$, which corresponds to the smallest power of two having 82 trailing non-\dig{2} digits; this same example has, in fact, 98 trailing non-\dig{2} digits. On the other hand, the total number of ternary digits of this power of two is about $1.3 \times 10^{14}$, far exceeding this 98 digit count.

\begin{figure}[t]
\centering
\includegraphics[scale=0.9]{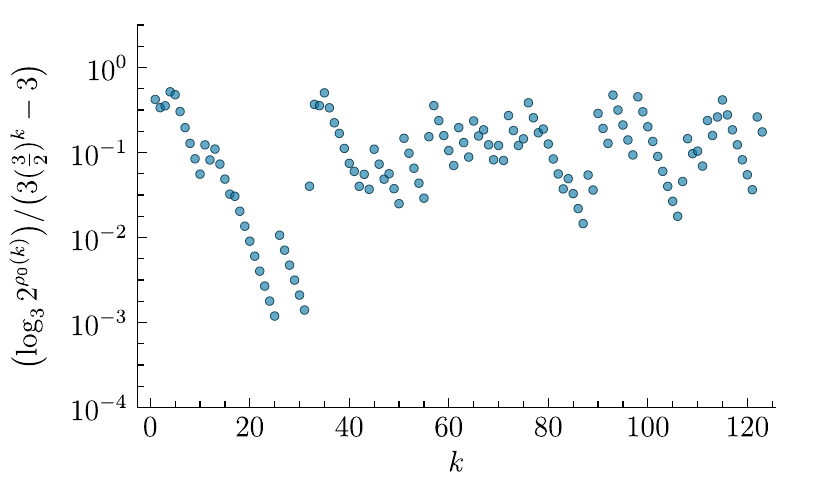}
\includegraphics[scale=0.9]{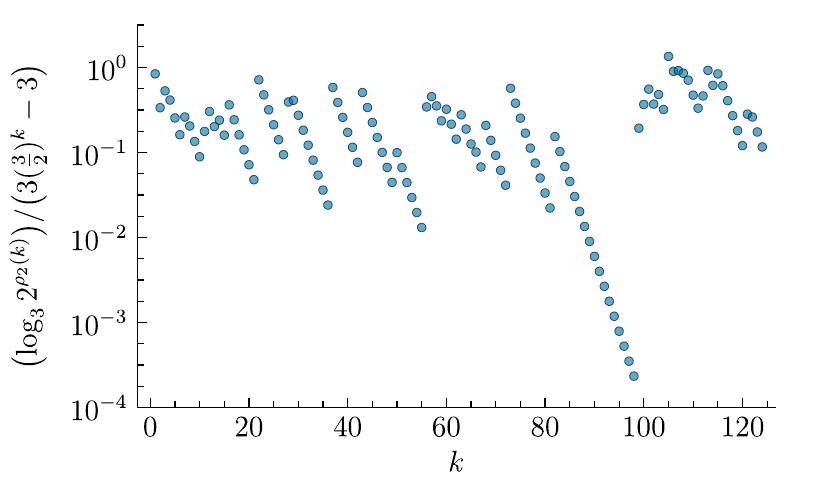}
\caption{Length of the ternary representation of $2^{\rho_\chi}$ (being approximately $\rho_\chi \log_3 2$) as a fraction of the expected number of rolls of three-sided die required to generate an uninterrupted sequence of $k$ non-$\chi$ digits (that average roll count being $3 ( \tfrac{3}{2} )^k - 3$).}
\label{fig2}\vspace{-1.5em}
\end{figure}

An alternative analysis comes from the heuristic that the ternary digits of powers of two are essentially random. Imagining the digits of $2^n$, reading from right-to-left, are a random number generator implementing the rolls of a three-sided die, we may ask how many rolls are necessary to generate an uninterrupted sequence of $k$ non-$\chi$ digits. Each non-$\chi$ digit has a probability of $\tfrac23$, and a routine calculation shows that we need, on average, $\smash{3 ( \tfrac{3}{2} )^k - 3}$ total rolls to generate such a sequence. Of course, this is only an approximation given that the digits of powers of two are entirely deterministic; in particular, the first and last digit of $2^n$ is never a \dig{0}, so this heuristic analysis could be slightly improved. Nevertheless, the expected roll count serves as an estimate of what the total ternary digit length is expected to be. Corresponding to the record breakers, \cref{fig2} plots the ternary digit length of $2^n$ (being approximately $n \log_3 2$) as a fraction of the expected roll count. We observe that, within zero to four of orders of magnitude, the digit counts of record breakers roughly match the expected roll count. The example of $n = 201015414581294$, mentioned in the previous paragraph, is uncharacteristic in the sense that for $k = 98$, we expect to require about $5.4 \times 10^{17}$ rolls, yet $2^{201015414581294}$ has only $1.3 \times 10^{14}$ ternary digits. Nevertheless, we observe in \cref{fig2} that there is no reasonable indication of finding any counterexamples to the conjectures: even the outlier record breakers are nowhere close to having the entire string of digits devoid of $\chi$.

\section{Conclusions}

By way of a recursive algorithm and extensive computation, we studied here the two conjectures of Erd\H{o}s and Sloane. These conjectures essentially state that, except for small number of trivial cases, every power of two has all possible digits somewhere in its ternary representation. The recursive algorithm explicitly constructs powers of two such that their trailing digits satisfy a certain requirement, e.g., consist solely of \dig0s and \dig1s. Testing these conjectures against all powers $2^n$ with $n \leq 2 \cdot 3^{45} \approx 5.9 \times 10^{21}$, no counterexamples were found. This extends an earlier study by Vardi \cite{vardi1991computational} which considered $n \leq 2 \cdot 3^{20} \approx 7 \times 10^9$. As part of the analysis, two ``record breaking'' integer sequences were defined: these record the smallest powers of two having no \dig{0} (resp., \dig{2}) in the last $k$ digits of its ternary representation, for $k = 1, 2, \ldots$. These integer sequences have been entered into the OEIS as \seqnum{A351927} and \seqnum{A351928}.

\appendix

\section{Proof of \texorpdfstring{\Cref{lemm}}{Lemma}}

We begin with a few elementary observations:
\begin{enumerate}[label=(\alph*)]
\item Suppose the last $k \geq 2$ ternary digits of an integer $x$ are $(a [\dig0]^{k-2} \dig1 )_3$ with $a \in \{\dig0,\dig1,\dig2\}$. (Here and in the following, the notation $[\,\cdot\,]^\ell$ means $\ell$ copies of the indicated digit.) Then, for some exponent $i \in \mathbb N$, we have that $x^i \equiv (a \cdot 3^{k - 1} + 1)^i \equiv ai \cdot 3^{k-1} + 1 \pmod{3^k}$, as shown by a simple application of the binomial theorem.

\item For a positive integer $k$, the last $k + 1$ ternary digits of $2^{u_k}$ are $(\dig1 [\dig0]^{k-1} \dig1)_3$. A simple inductive proof is as follows. Suppose the result holds for some $k \geq 2$ (the base cases with $k \in \{1,2\}$ are trivial to verify). Then $2^{u_k} - 1 = (3x + 1) 3^k$ for some non-negative integer $x$, and so
\begin{align*}
2^{u_{k+1}} &= (2^{u_k})^3 = \bigl((2^{u_k} - 1) + 1\bigr)^3 = (2^{u_k} - 1)^3 + 3 (2^{u_k} - 1)^2 + 3 (2^{u_k} - 1) + 1 \\
&\equiv 3^{k+1} + 1 \pmod {3^{k+2}},
\end{align*}
as required.

\end{enumerate}

\noindent Applying these observations, the proof of \cref{lemm} is as follows.

\begin{enumerate}
\item For $k \geq 2$, assume by induction that $u_{k-1}$ is the smallest positive integer $j$ such that $2^j \equiv_{k-1} 1$, and let $\ell$ be the smallest positive integer such that $2^\ell \equiv_k 1$. This number clearly satisfies $2^\ell \equiv_{k-1} 1$, and so if $\ell = a u_{k-1} + b$ with $a,b \in \mathbb N$ and $0 \leq b < u_{k-1}$, we see that $(2^{u_{k-1}})^a 2^b \equiv_{k-1} 1$. This linear congruence problem has a unique solution, namely $2^b \equiv_{k-1} 1$, which by the inductive hypothesis implies $b = 0$, and so $\ell$ is a multiple of $u_{k-1}$. By observation (b) above, $\ell$ cannot equal $u_{k-1}$ because the $k^\text{th}$ digit of $2^{u_{k-1}}$ is \dig{1}. Further, $\ell$ cannot equal $2 u_{k-1}$ because the square of $2^{u_{k-1}}$ has $k^\text{th}$ digit equal to \dig{2}. The next multiple of $u_{k-1}$ satisfies all requirements, and so $\ell = 3 u_{k-1} = u_k$, as claimed. (Note: the base cases of the inductive argument trivially hold by elementary computation.)

\item Suppose $i, j \in \mathbb N$ are such that $2^i \equiv_k 2^j$. Without loss of generality, suppose $i < j$. Then $2^i 2^{j - i} = 2^j$ yields a linear congruence $(2^i \bmod 3^k)(2^{j - i} \bmod 3^k) \equiv 2^j \pmod{3^k}$. Since the gcd of $(2^i \bmod 3^k)$ and $3^k$ is unity, there is exactly one solution to the linear congruence, namely $2^{j-i} \equiv_k 1$. Now suppose $j - i = a u_k + b$ with $a,b \in \mathbb N$ and $0 \leq b < u_k$; since $2^{j-i} = (2^{u_k})^a 2^b \equiv_k 1$ and $2^{u_k} \equiv_k 1$, again by uniqueness of the linear congruence problem, we find that $2^b \equiv_k 1$. Part (i) then implies $b = 0$, and so $i$ and $j$ differ by a multiple of $u_k$, as claimed.

\item Suppose $i, j \in \mathbb N$. Note that $2^{i u_k + j} \equiv_{k+1} (2^{u_k} \bmod 3^{k+1})^i 2^j \pmod{3^{k+1}}$. By observations (a) and (b), the trailing $k+1$ ternary digits of the first term are $\bigl([i \bmod 3][\dig0]^{k-1}\dig1\bigr)_3$. It is then a straightforward application of long multiplication to show that, modulo three, the $(k+1)^\text{st}$ digit of $2^{i u_k + j}$ is equal to the sum of the $(k+1)^\text{st}$ digit of $2^j$ plus $i$ times the first digit of $2^j$, as claimed.
\end{enumerate}

\section*{Acknowledgements}

The author thanks an anonymous reviewer for suggesting refinements to the proof of \cref{lemm}. Some computations used resources made possible by the Applied Mathematics Program of the U.S.~Department of Energy Office of Advanced Scientific Computing Research under contract number DE-AC02-05CH11231.

\bibliographystyle{siamplain}
\bibliography{references}

\end{document}